\documentclass{notices}
\usepackage{amssymb}
\usepackage{amsthm}
\usepackage[english]{babel}
\usepackage{tikz}
\usepackage{xcolor}
\colorlet{DarkGreen}{green!50!black}
\usepackage{float}
\usepackage{amsmath}

\newtheorem*{theorem}{Theorem}
\newtheorem*{conjecture}{Conjecture}
\newtheorem*{vgc}{Vague General Conjecture}

\newcommand{\SL}{\mbox{SL}}
\newcommand{\lieg}{\mathfrak{g}}
\newcommand{\lieh}{\mathfrak{h}}
\newcommand{\SO}{\mbox{SO}}
\newcommand{\Isom}{\mbox{Isom}}
\newcommand{\Conf}{\mbox{Conf}}

\newcommand{\Ad}{\mbox{Ad }}
\newcommand{\BR}{{\bf R}}
\newcommand{\BZ}{{\bf Z}}
\newcommand{\BC}{{\bf C}}

\title{
Rigidity of Transformation Groups in Differential Geometry
}

\author{
 Karin Melnick
  \affil{I thank David Fisher, Spyridon Lentas, and Vladimir Matveev
    for helpful feedback on a previous draft of this survey.}
 }

\begin{document}

\maketitle

\section*{Introduction and Framework}

A generic Riemannian manifold has no isometries at all.  But
mathematicians and physicists do not
really study such manifolds.  The first Riemannian manifolds we
learn about, and the most
important ones, are Euclidean space, hyperbolic space, and the sphere,
followed by quotients of these spaces, and by other symmetric
spaces.  These spaces go hand in hand with Lie groups and their
discrete subgroups.

In this survey, symmetry provides a framework for classification of
manifolds with differential-geometric structures.  We highlight
pseudo-Riemannian metrics, conformal structures,
and projective structures.  A range of techniques have been
developed and successfully deployed in this subject, some of them
based on algebra and dynamics and some based on analysis.  We aim to
illustrate this variety below.  



\subsection*{Historical Currents}

As with much of modern geometry, this subject has
intellectual foundations in the 1872 Erlangen Program of F. Klein.   He
proposed to study geometries via their transformation groups.  Euclidean
geometry thus comprises the objects and quantities invariant by translations, rotations, and
reflections, while conformal geometry is further enriched with dilations and circle
inversions.  Klein's classical geometries are homogeneous
spaces---their transformation groups act transitively---, and
the groups are Lie groups.

The concept of Lie group
arose from the
search for a theory of symmetries of differential equations.  In the
1930s, \'E. Cartan laid several foundations of differential geometry
with Lie groups: he found a complete set of local differential
invariants for many important and interesting geometric structures via
his method of \emph{moving frames}.
He defined the notion of \emph{$G$-structure of finite type} and his
beautiful theory of \emph{Cartan connections}, in which a manifold is
``infinitesimally modeled'' on one of Klein's homogeneous spaces.

Several decades later,
geometers such as S. Sternberg, V. Guillemin, and S. Kobayashi proved
that the local symmetries of any $G$-structure of finite type form a
finite-dimensional Lie pseudogroup.  This notion was then
greatly extended by Gromov in the 1980s with his \emph{rigid
  geometric structures}, which arose from his theory of partial
differential relations.

Meanwhile, in the study of discrete subgroups of Lie groups, A. Weil,
G. Mostow, and, ultimately, G. Margulis discovered remarkable rigidity
of lattices in semisimple Lie groups.  Among the breakthroughs in
these monumental results was the surprisingly powerful application of ergodic theory.  

Notions of rigid geometric structures will be explained and partially
defined below.  The reader could opt to read this section now, to read
it later, or to
skip it.  The article can be well appreciated by just
considering the examples of geometric structures given. 

\subsection*{Zimmer's Program}
In
the 1970s and 80s, R. Zimmer developed rigidity theory of Lie groups and their
lattices in the context of their actions on manifolds (see \cite{zimmer.etsg}).  A
\emph{lattice} in a Lie group $H$ is a discrete subgroup $\Gamma$ such
that the quotient $H / \Gamma$ has finite Haar measure.

Let $\Gamma$ be a lattice in a simple Lie group $H$; we assume
for this subsection and the next that $H$ is noncompact, has finite center, and is not locally
isomorphic to $\mbox{SO}(1,n)$ or $SU(1,n)$, meaning the Lie algebra
is not isomorphic to that of either of these groups.  Superrigidity in this
case says that homomorphisms of $\Gamma$ to
linear Lie groups $G$ are the restrictions of homomorphisms $H
\rightarrow G$, up to some ``precompact noise,'' and possibly
passing to a cover of $H$.  The linear manifestations of these
complicated and mysterious objects are thus reduced to those of
$H$---an elementary algebraic matter, solved by Schur's Lemma.
(Superrigidity and the remaining results and conjectures
of this section are valid for $H$ semisimple, but this entails additional definitions and assumptions.)

Zimmer's Program asserts that, in many cases, homomorphisms from $\Gamma$, or $H$, as above to
the group of volume-preserving diffeomorphisms of a compact manifold
arise from a short list of algebraic constructions.   We illustrate
two such actions, along with invariant
differential-geometric structures for each.  Let $M = G/\Lambda$, for $G$ a connected Lie group and $\Lambda < G$
  a cocompact lattice. 

\begin{enumerate}
\item A homomorphism $\Gamma \rightarrow G$ gives
  an action of $\Gamma$ on $M$, preserving the volume determined by the Haar measure.

  For $G$ semisimple, the \emph{Killing form} on the Lie algebra
$\lieg$ is a nondegenerate bilinear form, invariant by conjugation via the adjoint representation, $\Ad G$. It gives rise to a
bi-invariant metric on $G$, of indefinite signature when $G$ is
noncompact.  
This \emph{Cartan-Killing metric}
 descends to a $\Gamma$-invariant pseudo-Riemannian metric on $M$.

  \item Any subgroup $\Gamma \leq \mbox{Aut}_{\Lambda} G$, the automorphisms of
    $G$ normalizing $\Lambda$, acts on $M$.  The standard example is
    $G/\Lambda = \BR^n / \BZ^n = {\bf T}^n$, and $\Gamma =
    \mbox{SL}_n(\BZ)$.

    Any Lie group $G$ has a bi-invariant affine connection $\nabla_G$, defined by
    declaring left-invariant vector fields to be parallel.  This
    connection descends to $M$ and is $\Gamma$-invariant.  Although
    $\nabla_G$ is torsion-free only if $G$ is abelian, the existence
    of a $\Gamma$-invariant affine connection on $M$ implies existence
    of a torsion-free affine connection that is also
    $\Gamma$-invariant.
  \end{enumerate}


  Pseudo-Riemannian metrics and affine connections are examples of $G$-structures of finite type.
For $\Gamma$ a lattice in a simple Lie group $H$ as above, an
\emph{affine action} of $\Gamma$ is an action on $G / \Lambda$ via a
homomorphism to $\mbox{Aut}_{\Lambda} G \ltimes G$.  We can refine the
classification claim from  Zimmer's Program to say,
\emph{infinite actions of $\Gamma$ or $H$ on a compact manifold $M$ preserving a
$G$-structure of finite type and a volume are affine actions---up to some additional ``precompact
noise,''} for which we refer to \cite{fisher.survey.festschrift} for details.

A famous conjecture of Zimmer is a dimension bound, with no invariant
rigid geometric structure assumed:
\emph{Let $N$ be the
minimal dimension of a nontrivial representation of $H$.   For any
compact manifold $M$ of dimension $n < N$, any action of $\Gamma$ on
$M$ by volume-preserving diffeomorphisms preserves a smooth Riemannian
metric.}  Depending on $H$, the conclusion can often be strengthened
to the action being finite.
Affine actions
satisfy this conjecture, which can be proved with Zimmer's
Cocycle Superrigidity Theorem.  
Given a volume-preserving action of $\Gamma$ on any $n$-dimensional
manifold, for any $n$, Cocycle Superrigidity gives a nontrivial homomorphism $H \rightarrow
\mbox{SL}_n(\BR)$ or a measurable, $\Gamma$-invariant
Riemannian metric on $M$.  In the case $\Gamma$ also preserves
a $G$-structure of finite type, Zimmer improved the measurable
Riemannian metric to a smooth one.

For decades, researchers struggled to improve the invariant, measurable
Riemannian metric in the general case.  The
breakthrough came in 2016, when Brown--Fisher--Hurtado
proved the
dimension conjecture for $\Gamma$ cocompact and $H$ any classical
simple Lie group, aside from a few exceptional dimensions for types
$B_n$ and $D_n$.
They also proved a dimension bound for smooth actions not necessarily
preserving a volume: in this case, $\mbox{dim } M \geq N-1$, or the
action is finite.  A sharp example to have in mind is the action of
$\mbox{SL}_n(\BR)$, or one of its lattices, on $M = {\bf RP}^{n-1}$.
Their proof opens up new avenues for the classification conjecture of
Zimmer's Program---see below.

\subsection*{Zimmer-Gromov Program}

 Gromov generalized the $G$-structures of finite type in
Zimmer's earlier work to the richer \emph{rigid geometric
  structures}.  Inspired in large part by Zimmer's work, he proved
results for a wider class of Lie groups with actions preserving
these structures.  His vision to some extent parts ways with another, large branch of Zimmer's
Program, focused on
rigidity of lattices in semisimple Lie groups, which
has been demonstrated for a breathtaking range of actions.

Gromov attempted to formulate Zimmer's conjecture about affine
actions more broadly, in particular, without assuming a finite,
invariant volume.
In \cite{dag.rgs}, he and D'Ambra state the
\begin{vgc}[VGC]
All triples $(H,M,\omega)$, where $M$
is a compact manifold with rigid geometric structure $\omega$ and H is a
``sufficiently large'' transformation group, are almost classifiable.
\end{vgc}
\noindent
``Large'' here often means noncompact.  It could also be
a transitivity condition, such as having a dense orbit or an ergodic
invariant measure.  It could be a geometric condition of being \emph{essential},
generally meaning it does not preserve a
finer geometric structure subsidiary to $\omega$. 
They note, ``we are still far from proving (or even
starting) this conjecture, but there are many concrete results which
confirm it.''  Some such results will be presented below.

Here are examples of actions that join the
``almost classification'' for large actions of arbitrary
groups not necessarily preserving a volume.  For a more complete list, see \cite[0.5]{gromov.rgs}

\begin{enumerate}

\item[3.] $M = G/P$ with $P$ closed---that is, an arbitrary
  homogeneous space.  Assuming $G$ semisimple, an important class are the
  quotients by \emph{parabolic subgroups}, algebraic subgroups of $G$
  for which $G/P$ is a compact projective variety.  These $G$-actions
  never preserve a finite volume.
  
The basic example is projective space ${\bf RP}^n$, with $G =
\mbox{PSL}_{n+1}(\BR)$, and $P$ the stabilizer of a line in
$\BR^{n+1}$.  This is a Klein geometry, in which the invariant notions are
lines and intersection.

Another important parabolic Klein geometry is the round sphere
${\bf S}^n$, a homogeneous space
of $\mbox{PO}(1,n+1) \cong \mbox{Isom}({\bf H}^{n+1})$.  It is the
projectivization of the light cone in Minkowski space
${\BR}^{1,n+1}$.  It carries an invariant conformal structure, for which the invariant notion is
angle. 
  
  \item[4.] $M = \Gamma \backslash G/P$, where $\Gamma < G$ acts freely, properly
    discontinuously, and cocompactly on $G/P$.

    An example is the unit
    tangent bundle of a compact hyperbolic manifold, where
    $G \cong \mbox{PO}(1,n+1)$ and $P \cong \mbox{O}(n)$ is the stabilizer of a unit
    vector in $T{\bf H}^{n+1}$.  The image of the holonomy
    representation of $\pi_1(M)$ is $\Gamma$.  The geodesic flow on
    $M$ is furnished by a one-parameter subgroup of $G$ centralizing
    $P$; it is isometric for a metric obtained by pseudo-Riemannian
    submersion from the
    Cartan-Killing metric on $\Gamma \backslash G$.

\item[5.] $M = \Gamma \backslash U$, where $U \subset G/P$ is open, and $\Gamma$ acts freely, properly discontinuously, and
  cocompactly on $U$.  The example of the conformal Lorentzian Hopf
  manifold is presented later.

The spaces in (4) and (5) both carry \emph{$(G,{\bf X})$-structures}, which are defined in the next section.
  \end{enumerate}

  Gromov's first key contribution, based on his theory of
  partial differential relations, is known as the Frobenius Theorem.
 At points satisfying a regularity condition, it can produce local transformations of $(M,\omega)$ from pointwise, infinitesimal input.  The infinitesimal
  input corresponds to points in an algebraic variety.  The Rosenlicht
  Stratification for algebraic actions on varieties yields a corresponding
  stratification by local transformation orbits in $M$.
  A consequence is the Open-Dense Theorem: if an orbit for local
  transformations in $(M,\omega)$ is dense,  
then an open, dense subset $U \subseteq M$ has a $(G,{\bf X})$-structure.

The Open-Dense Theorem lends support to the VGC.  In specific cases,
researchers can show that
 $U = M$, and some have obtained
significant theorems in the Zimmer-Gromov Program by this route.  But $U$ may not,
in general, equal $M$.  Classifying compact $(G,{\bf
  X})$-manifolds compatible with a given geometric structure is a
challenging subject unto itself.

A beautiful alternative relating examples 1-2 with example 3 is given by results of Nevo--Zimmer on smooth projective
factors.  They use Gromov's Stratification and his Representation
Theorem to prove: given a connected, simple 
$H$ preserving a real-analytic rigid geometric structure
on a compact manifold $M$, together with a
\emph{stationary measure}, there is $i$) a smooth, $H$-equivariant projection
of an
open, dense $U \subseteq M$ to a parabolic homogeneous space $H/Q$; or
$ii$) the Gromov representation of $\pi_1(M)$ contains $\lieh$ in its
Zariski closure.  A stationary measure is a finite measure, here
assumed of full support, invariant by
convolution with a certain finite measure on $H$, but not necessarily
under the $H$-action.

\subsection*{Notions of rigid geometric structure}

The infinitesimal data of a geometric structure on $M$ inhabits the frame
bundle of $M$, of the appropriate order.  The frames at a point $x \in
M$ are the linear isomorphisms $T_xM \rightarrow \BR^n$, where $n =
\mbox{dim } M$.  These form a principal $\mbox{GL}_n(\BR)$-bundle,
denoted $\mathcal{F}M$.
For $G \leq \mbox{GL}_n(\BR)$, a $G$-structure is a \emph{reduction} $\mathcal{R}$
of $\mathcal{F}M$ to $G$---that is, a smooth $G$-prinicipal
subbundle.  The frames in $\mathcal{R}$ determine the
structure.  For example, a $(p,q)$-semi-Riemannian metric $g$ is a
$G$-structure for $G \cong \mbox{O}(p,q)$.  The frames of $\mathcal{R}$ are
the bases
of $T_xM$ in which the inner product $g_x$ has the form
$-\mbox{Id}_p \oplus \mbox{Id}_q$.

The fiber $\mathcal{F}_xM$ comprises the first derivatives of coordinate
charts at $x$.  The order-$k$ frame bundle has fibers
$\mathcal{F}^{(k)}_xM$ comprising $k$-jets of coordinate charts $\varphi$
at $x$ with $\varphi(x) = {\bf 0}$.
Here $\varphi$ and $\psi$ have the
same $k$-jet at $x$ if $D^{(i)}_{\bf 0}(\varphi \circ \psi^{-1}) =
D^{(i)}_{\bf 0}
\mbox{Id}$ for all $1 \leq i \leq k$.  The $k$-jets at ${\bf 0}$ of
diffeomorphisms of $\BR^n$ fixing ${\bf 0}$ form the group
$\mbox{GL}^{(k)}_n(\BR)$.  A $G$-structure of order $k$ is a reduction
of $\mathcal{F}^{(k)} M$ to $G \leq \mbox{GL}_n^{(k)}(\BR)$.  For
example, a linear connection on $TM$ is a reduction
$\mathcal{R} \subset \mathcal{F}^{(2)} M$ to $\mbox{GL}_n(\BR) < \mbox{GL}_n^{(2)}(\BR)$.  
Elements of $\mathcal{R}$ are $2$-jets of geodesic coordinate charts.

A $G$-structure is \emph{finite type} if a certain prolongation
process stabilizes, which essentially means there is $k$ such
that, at all $x \in M$, the $k+i$-frames at $x$ adapted to the
structure are determined by the adapted
$k$-frames at $x$, for all $i \geq 0$ (see \cite{kobayashi.transf}).  Conformal Riemannian structures, for example,
are $\mbox{CO}(n)$-reductions of $\mathcal{F}M$, where $\mbox{CO}(n) \cong
\BR^* \times \mbox{SO}(n)$.  The prolongation gives a
reduction of $\mathcal{F}^{(2)}M$ to a bigger group, namely the
parabolic subgroup $P$ fixing a point of ${\bf S}^n$.  \'E. Cartan
found, for $n \geq 3$, the prolongation stabilizes here, and thus any conformal
transformation is determined by the 2-jet at a point.

A $G$-reduction of $\mathcal{F}M$ is the same as an
equivariant map $\mathcal{F}M \rightarrow \mbox{GL}_n(\BR)/G$ (and similarly in higher order).
Gromov's \emph{geometric structures of algebraic type} are equivariant maps
to algebraic varieties.  They are \emph{rigid} if a similar prolongation
process stabilizes (see \cite{gromov.rgs}).  Some very useful additional flexibility is
afforded by varieties which are not necessarily
homogeneous.  For example, one may add to a $G$-structure of finite
type some vector fields, which correspond to maps to $\BR^n$, and the
result is again a rigid geometric structure.  All rigid
geometric structures in this article are understood to be of algebraic type.

For $G$ a Lie group and $P <
G$ a closed subgroup, a \emph{Cartan geometry} on $M$
modeled on $(G,P)$ comprises a principal $P$-bundle $
\pi : \hat{M} \rightarrow M$ and a $\lieg$-valued $1$-form $\omega$ on
$\hat{M}$, called the \emph{Cartan connection}.  The Lie group $G$ carries a left-invariant
$\lieg$-valued $1$-form, the \emph{Maurer-Cartan form} $\omega_G$,
which simply identifies the left-invariant vector fields on $G$ with
$\lieg$.  The Cartan connection is required to satisfy three axioms,
mimicking properties of $\omega_G$.  One of them says $\omega$ is an
isomorphism on each $T_{\hat{x}} \hat{M}$, so it
determines a parallelization of $T\hat{M}$ (see \cite{sharpe}).  

Essentially all classical rigid geometric structures canonically
determine a Cartan geometry, such that the transformations correspond
to automorphisms of the Cartan geometry---diffeomorphisms of $M$
lifting to bundle automorphisms of $\hat{M}$ preserving $\omega$.  A
conformal Riemannian structure in dimension $n \geq 3$, for example, determines a Cartan
geometry modeled on ${\bf S}^n$, with $G = \mbox{PO}(1,n+1)$ and $P$ as
above.  The bundle $\hat{M}$ in this case is the $P$-reduction of
$\mathcal{F}^{(2)}M$ obtained by prolongation
of the $\mbox{CO}(n)$-structure.

Cartan geometries have proven very useful for the study of
transformation groups.  Fundamental results of the Zimmer-Gromov theory of
rigid geometric structures have been translated to this setting.  A key asset is the
\emph{curvature 2-form}
$$ \Omega = d \omega + \frac{1}{2} [ \omega, \omega ] \in
\Omega^2(\hat{M}, \lieg)$$
Vanishing of $\Omega$ over an open $U \subseteq M$ is the
obstruction to $U$ having a $(G,{\bf X})$-structure with
${\bf X} = G/P$.

The notion of \emph{$(G,{\bf X})$-structure}
   was developed by Ehresmann and later by Thurston \cite{thurston.3d.book}: on a manifold $M$, it comprises an atlas of charts to ${\bf X}$
    with transitions equal to transformations in $G$ (one $g \in G$ on
    each connected component of the chart overlap).  When $M = \Gamma
    \backslash {\bf X}$ for some $\Gamma < G$, it is called
    \emph{complete}; quotients $\Gamma \backslash U$ of open subsets
    $U \subset {\bf X}$ are called \emph{Kleinian}.

\section*{Isometries of pseudo-Riemann-ian manifolds}
By the classical theorems of Myers and Steenrod, the isometry group of
a compact Riemannian manifold is a compact Lie group.  Isometries of
pseudo-Riemannian manifolds need
not act properly; in particular, when the manifold is compact, this group can be noncompact.
In this section we describe some important results on isometry groups
and illustrate  a couple techniques which have been influential.
We focus on compact Lorentzian manifolds, mostly
because there are as yet few answers to the corresponding questions in
higher signature.  The term \emph{semi-Riemannian} means Riemannian or pseudo-Riemannian.

\subsection*{Simple group actions: Gauss maps}

Let $G$ be a noncompact, simple Lie group.  For $\Lambda < G$ a lattice, the
$G$-action on $G/\Lambda$ is
\emph{locally free}, meaning stabilizers are discrete.
We sketch an argument of Zimmer, using
the Borel Density Theorem, that 
any isometric $G$-action on a finite-volume pseudo-Riemannian
manifold $M$
is locally
free on an open, dense subset.
First suppose that $G$ acts
ergodically, meaning
that any $G$-invariant measurable function is constant
almost everywhere.  In fact, any $G$-invariant measurable map to a
\emph{countably separated} Borel space is constant almost everywhere;
countably separated means there is a countable collection of Borel
subsets $B_i$ such that for any points $x \neq y$, some $B_i$ contains
one of $\{ x,y\}$ and not the other.  This property holds for quotients of algebraic varieties by algebraic actions.
Then we can apply this fact to
$G$-equivariant maps from
$M$ to $G$-algebraic varieties, sometimes called \emph{Gauss maps}.

Let $\sigma$ assign to $x \in M$ the Lie algebra of the stabilizer of
$x$ in $G$, which lies in the disjoint union $\mbox{Gr } \lieg$ of the
Grassmannians $\mbox{Gr}(k,\lieg)$, $k=0, \ldots, \mbox{dim }
\lieg$.
The map $\sigma$ satisfies the equivariance relation $\sigma(g.x) = (\mbox{Ad } g). (\sigma(x))$.  We conclude
that there is a single $G$-orbit $\mathcal{O} \subset \mbox{Gr }
\lieg$ containing $\sigma(x)$ for almost-every $x$.
The volume on $M$ pushes forward to an $\Ad G$-invariant finite volume
on $\mathcal{O}$. After passing to Zariski closures, this
orbit can be identified with an algebraic homogeneous space.
 The Borel Density Theorem says that $\overline{\mathcal{O}}$, hence
 also the $\Ad G$-orbit $\mathcal{O}$, must be
 a single point.  Thus $\sigma(x)$ is an ideal
for almost-all $x \in M$, namely, because $G$
is simple, it is $0$ or $\lieg$.

In general, the metric 
volume decomposes into ergodic components, and $\sigma(x) \in \{ 0,
\lieg \}$ for almost-all $x$.  The points with $\sigma(x) = \lieg$ comprise the fixed set
of $G^0$.  The fixed set of any nontrivial isometry has null volume
(assuming $M$ connected).  Thus
$\sigma(x) = 0$, and the stabilizer of $x$ is discrete, for all $x$ in
a full-volume subset $\Omega \subset M$.  It follows from lower semicontinuity of
$\mbox{dim } \sigma(x)$ that $\Omega$ is open and dense.

By comparable arguments involving ergodicity and the Borel Density
Theorem, Zimmer obtained, for $M$ compact, a Lie algebra embedding $\lieg
\hookrightarrow \mathfrak{o}(p,q)$ and concluded that, in particular,
the real-rank---the dimension of a maximal abelian,
$\BR$-diagonalizable subgroup of $\Ad G$---satisfies
$$\mbox{rk}_\BR G \leq \mbox{min} \{p,q\} =\mbox{rk}_{\BR} \mbox{SO}(p,q)$$
Such an embedding for $G$ of real-rank at least two
follows from Zimmer's Cocycle Superrigidity; here, with a
pseudo-Riemannian metric, a stronger result is obtained
by a more elementary proof.

A more general embedding theorem for connected $G$, not
necessarily simple, preserving any rigid geometric structure, not necessarily determining a volume, was proved by
Gromov \cite{gromov.rgs}.  A version of this generality was proved for
compact Cartan geometries by Bader--Frances--Melnick in 2009. 

\subsection*{Connected isometry groups of compact Lorentzian manifolds}

The group $\SL_2(\BR)$ with the Cartan-Killing metric is a
three-dimensional Lorentzian manifold with isometry group $\SL_2(\BR)
\times_{\BZ_2} \SL_2(\BR) \cong \SO^0(2,2)$.  This space has constant
negative sectional curvature and is known as anti-de Sitter space, $\mbox{AdS}^3$. (Together with
its higher-dimensional analogues, this space is important in string theory and M-theory.)

Let $(M,g)$ be a compact, connected, Lorentzian manifold and $H \leq
\Isom(M,g)$ a noncompact simple group.  Zimmer's bound is
$\mbox{rk}_\BR H \leq 1$; he in fact proved that the identity
component $\Isom^0(M,g)$
is
locally isomorphic to $\SL_2(\BR) \times K$, with
$K$ compact.  Gromov improved Zimmer's result to conclude
$$M \cong (\widetilde{\mbox{AdS}}^3 \times_f N)/ \Gamma$$
where $N$ is a Riemannian manifold, $f : N \rightarrow \BR^+$ is a
warping function, and $\Gamma$ acts isometrically, freely, and
properly discontinuously on the warped product.

The complete determination of connected isometry groups of compact
Lorentzian manifolds was simultaneously achieved by Adams--Stuck
\cite{as.lorisom1, as.lorisom2} and Zeghib
\cite{zeghib.lorisom1, zeghib.lorisom2}:

\begin{theorem}[Adams--Stuck/Zeghib 1997/8]
The identity component of the isometry group of a compact Lorentzian
manifold is locally isomorphic to $H \times K
\times \BR^k$, with $K$ compact, $k \geq 0$, and $H$ one of
\begin{itemize}
\item $\SL_2(\BR)$
\item $\mbox{Heis}^m$, a $2m+1$-dimensional Heisenberg group
\item an oscillator group, a solvable extension of
  $\mbox{Heis}^{2m+1}$ by $S^1$.
\end{itemize} 
\end{theorem}

Adams and Stuck proceed by studying the dynanics of a Gauss map
$M \rightarrow \mbox{Sym}^2 \lieg$, where $\lieg$ is the Lie algebra of
$\Isom^0 M$.  The map sends $x \in M$ to the pullback of the metric on
$T_xM$ restricted to the subspace tangent to the $G$-orbit at $x$.  It
is equivariant for the $\Ad G$ representation on $\mbox{Sym}^2
\lieg$.  Zeghib works with the average over $M$ of these forms on $\lieg$.

The geometry of $M$ when $H$ is a Heisenberg or oscillator group has
not been completely described.

\subsection*{Dynamical foliations from unbounded isometries}

The anti-de Sitter space $\mbox{AdS}^3$ can be identified with the unit
tangent bundle of the hyperbolic plane.  For a hyperbolic
surface $\Sigma \cong \Gamma \backslash {\bf H}^2$, the quotient
$\Gamma \backslash \mbox{AdS}^3$ is identified with $T^1\Sigma$, and the geodesic and horocycle flows
with the right-action of the diagonal and upper-triangular unipotent
subgroups, respectively, of $\SL_2(\BR)$.

\begin{figure}[h]
  \begin{center}
    \includegraphics[width=3cm]{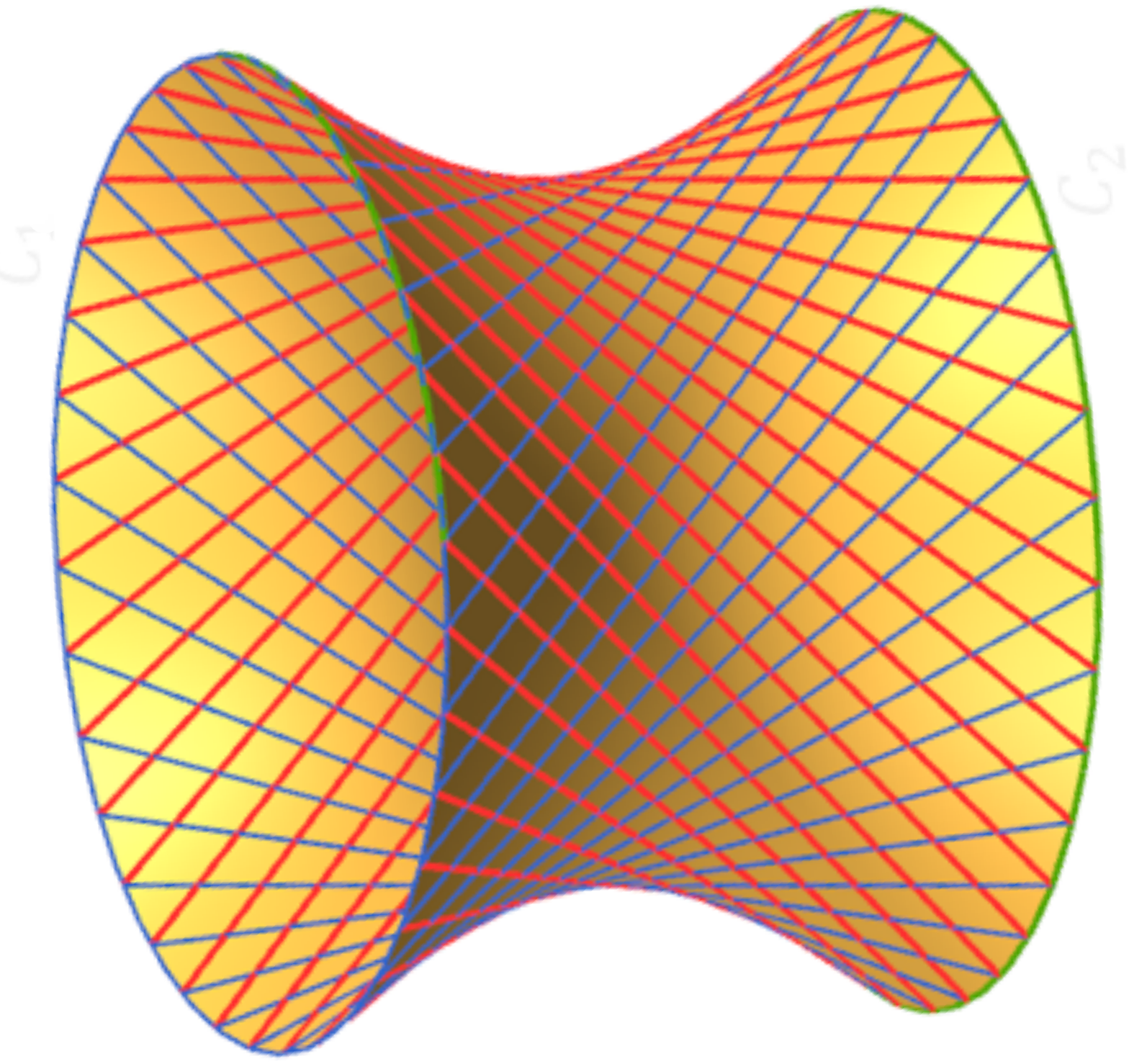}
    \caption{\label{fig.AdS.biruling} Two-dimensional Anti-de Sitter
      space with bifoliation by lightlike geodesics (creative commons:
      bubba on \url{math.stackexchange.com/q/1606820})}
  \end{center}
  \end{figure}

The geodesic flow is an Anosov, Lorentz-isometric flow on
$\Gamma \backslash \mbox{AdS}^3$.  The weak-stable and -unstable
foliations are by
totally geodesic surfaces on which the metric is degenerate.  These foliations, $\mathcal{W}^{cs}$ and $\mathcal{W}^{cu}$, are the orbits of the two-dimensional upper-triangular
and lower-triangular subgroups, respectively (see Figure \ref{fig.AdS.biruling}).  The horocycle flow is
not hyperbolic
in the dynamical sense.  Nonetheless, it also preserves a
foliation by totally geodesic, degenerate surfaces, namely $\mathcal{W}^{cs}$.  Zeghib calls this
the \emph{approximately stable} foliation of the flow.  He proved the
following \cite{zeghib.tgl1, zeghib.tgl2}:

\begin{theorem}[Zeghib 1999]
Let $M$ be a compact Lorentzian manifold, and suppose that
$G = \mbox{Isom }M$ is noncompact.  Any unbounded sequence in $G$
determines, after passing to a subsequence, a foliation by totally geodesic, degenerate hypersurfaces.
\end{theorem}
\noindent
The space of all these approximately stable foliations gives a compactification of $G$
and reflects algebraic properties of $G$.

The construction of such foliations was suggested by Gromov, who
provided the following argument.  The reader may note the significance
of the index $p=1$, a way in which Lorentzian geometry is ``closest'' to
Riemannian. Let $f_k \rightarrow \infty$ in $G$.  For any $x \in M$, we may
pass to a subsequence so that $f_k.x$ converges to a point $y \in M$.
Equip $M \times M$ with the pseudo-Riemannian metric $g \oplus -g$.
The Levi-Civita connection on the product is $\nabla^g \oplus
\nabla^g$.
The graphs $\Gamma_k$ of $f_k$ are totally isotropic and totally geodesic in
$(M \times M, g \oplus -g)$, and they converge to an $n$-dimensional,
isotropic, totally geodesic submanifold $\Gamma$ near $(x,y)$.

Because $\{f_k\}$ does not converge, $\Gamma$ cannot be a graph near $(x,y)$, so it
has positive-dimensional intersection with $\{ x \} \times M$.  This
intersection is isotropic and geodesic in $(M,g)$, so coincides with
an isotropic geodesic segment $\gamma$ near $y$.  The projection of $T_{(x,y)}
\Gamma$ to $T_y M$ must be orthogonal to $T_y \gamma$, so it is
contained in the degenerate hyperplane $T_y \gamma^\perp$.
Therefore, $T_{(x,y)} \Gamma$ has nonzero intersection with $T_x M \oplus
0$, which must again must be isotropic; the projection to $T_xM$
is contained in a degenerate hyperplane $H$.  A dimension count
implies that the
projection equals $H$.  Because $\Gamma$
is totally geodesic, it projects to the totally geodesic, degenerate
hypersurface obtained by exponentiation of $H$.

A classical theorem of Haefliger says that a compact, simply
connected, real-analytic manifold admits no real-analytic,
codimension-one foliation.  Zeghib used his foliations to give
another proof of the following theorem of D'Ambra \cite{dambra.lorisom}:

\begin{theorem}(D'Ambra 1988)
Let $n \geq 3$.  For $(M^n,g)$ a compact, simply connected, real-analytic Lorentzian
manifold, $\mbox{Isom}(M,g)$ is compact.
  \end{theorem}
  \noindent
  D'Ambra's proof was a \emph{tour de force} of Gromov's theory of rigid
geometric structures and belongs to the underpinning of the VGC.

\subsection*{Current Questions}

After reading about these wonderful accomplishments of the 1980s and
90s, the reader is likely to have at least two obvious questions: what
about nonconnected groups? and, what
about higher signature?  

Some progress has been made on compact Lorentzian manifolds $(M,g)$ for which
$G = \mbox{Isom}(M,g)$ has infinitely-many components when $M$ is
\emph{stationary}---that is, it admits a timelike Killing vector
field.  The latter property implies that $M$ has closed timelike geodesics,
or ``time machines.''
Piccione--Zeghib have a nice structure theorem for such spaces:
$G^0$ contains a torus ${\bf T}^d$, and there is a
Lorentzian quadratic form $Q$ on $\BR^d$ such that $G/G^0 \leq O(Q)_\BZ$.

For smooth, three-dimensional, compact Lorentzian manifolds, Frances
recently completely classified those with noncompact isometry group,
topologically and geometrically.  The topological classification says
that such a manifold $M$ is, up to a cover of order at most
four:

\begin{itemize}
\item $\Gamma \backslash \widetilde{\SL}_2(\BR)$, for $\Gamma$ a
  cocompact lattice
  \item The mapping torus ${\bf T}^3_A$ of an automorphism $A$ of
    ${\bf T}^2$.
\end{itemize}
In the first case, $\widetilde{\SL}_2(\BR)$ has a left-invariant
metric.  Geometries of the second type can have infinite, discrete
isometry group.  A corollary is an improvement of D'Ambra's Theorem to
smooth metrics, in the three-dimensional case.

 For $p \geq 2$, it is not known which
\emph{homogeneous}, compact, $(p,q)$-pseudo-Riemannian manifolds have noncompact
isometry group.  Much current work focuses on classifying left-invariant
metrics on Lie groups.  To illustrate the complexity of higher
signature, we indicate Kath--Olbrich's survey on the classification problem for
pseudo-Riemannian symmetric spaces \cite{kath.olbrich.survey.esi}.  While Riemannian symmetric spaces
were essentially classified by \'E. Cartan, and reductive
pseudo-Riemannian symmetric spaces by Berger in 1957, a complete classification
is not really expected, without restricting to
a low index $p$, or assuming some additional invariant structure.

\section*{Conformal Transformations}

Compact Lorentzian manifolds are marvelous for the mathematical
study of transformation groups, but the possibility of ``time
machines'' makes them less appealing to physicists.  Conformal
compactifications of Lorentzian manifolds are, on the other hand,
quite natural in relativity.

We have seen above that conformal Riemannian transformations are those
preserving angles between tangent vectors.  For arbitrary signature,
we define
the \emph{conformal class} of $g$ to be
$ [g] = \{ e^{2 \lambda} g \}$, with $\lambda$ ranging over smooth functions
on $M$.  The conformal transformations $\mbox{Conf}(M,[g])$ are those
preserving $[g]$.  In higher signature,
these are the transformations preserving the causal structure---that
is, the null cones in each tangent space.

Stereographic projection from a point $\hat{p} \in {\bf S}^n$ is a
conformal equivalence of ${\bf S}^n \backslash \{ \hat{p} \}$ with
Euclidean space.  Thus ${\bf S}^n$ is \emph{conformally flat}---locally conformally equivalent to
a flat Riemannian manifold.
The higher signature analogue of ${\bf S}^n$ is a
parabolic homogeneous space for $\mbox{PO}(p+1,q+1)$,
called the \emph{M\"obius space} ${\bf S}^{p,q}$.  It is a two-fold
quotient of $S^p \times S^q$; the conformal class can be realized by
$g_p \oplus -g_{q}$, where $g_p, g_q$ are the constant-curvature
metrics on the respective spheres.  These metrics
are also conformally flat.  The conformal structure of a
$(p,q)$-pseudo-Riemannian metric, $p+q \geq 3$, determines a Cartan
geometry modeled on ${\bf S}^{p,q}$.  It is also a $G$-structure of
finite type, with order
2.  Throughout this section, $n$ is an integer greater than or equal $3$.

\subsection*{The Ferrand-Obata Theorem}

\begin{figure}
  \begin{center}
    \includegraphics[width=.9in]{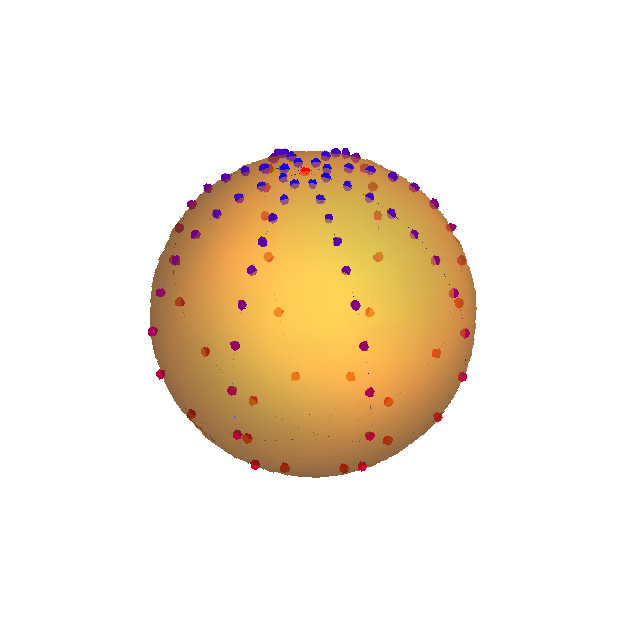} \ \ 
    \includegraphics[width=.9in]{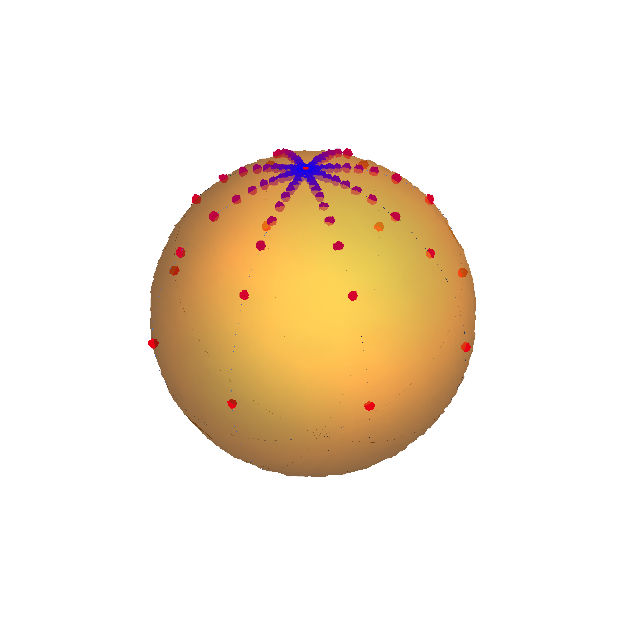} \ \ 
    \includegraphics[width=.9in]{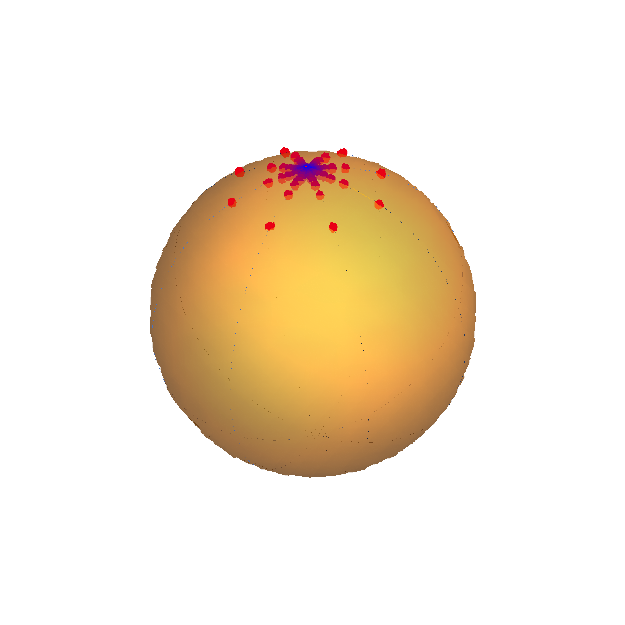}
  \caption{\label{figure.sphere} iterates of a conformal
    transformation of the
  sphere with source-sink dynamics}
  \end{center}
  \end{figure}
    
While the isometry group of a compact Riemannian manifold is compact,
the conformal group need not be.
A dilation of Euclidean space corresponds via stereographic projection
to a conformal transformation of ${\bf S}^n$ having
source-sink dynamics under iteration---see Figure
\ref{figure.sphere}.  A sequence $\{ f_k\}$ of transformations has
\emph{source-sink dynamics} if there are points $p^+$ and $p^-$ in $M$ such
that  $\{f_k\}$ converges uniformly
  on compact subsets of $M \backslash \{ p^-\}$ to the constant map $p^+$, 
  and similarly $f_k^{-1}\rightarrow p^-$ uniformly on compact subsets
 of $M \backslash \{ p^+ \}$.
The full conformal group of  ${\bf S}^n$ is the noncompact simple Lie
group $\mbox{Isom } {\bf
  H}^{n+1} \cong \mbox{PO}(1,n+1)$.
  This isomorphism comes from the
fact that ${\bf S}^n$ is the visual
boundary in the conformal compactification of ${\bf H}^{n+1}$.
Being the visual boundary of hyperbolic space is quite a special
property.  Lichnerowicz conjectured the following striking fact:

\begin{theorem}[Ferrand/Obata 1971]
Let $(M^n,g)$ be a compact Riemannian manifold.  If $\Conf(M,[g])$ is noncompact, then $(M,[g])$ is globally conformally
equivalent to ${\bf S}^n$.  
  \end{theorem}
  \noindent
  (This theorem is also true for $n=2$.)
Lichnerowicz's original conjecture was actually under the stronger
assumption $\Conf^0(M,[g])$ noncompact, and this was proved by Obata \cite{obata.lich}, using Lie
theory and differential geometry.  Ferrand proved the theorem in full
generality \cite{lf.lich}, using quasiconformal analysis and no Lie theory.

In 1994, Ferrand extended her result to noncompact $M$ and gave a
streamlined proof of her original result.
She uses
conformal capacities to define a conformally invariant function $\mu$ on
pairs of points.
Sometimes it
gives a metric on $M$; in this case, $\Conf(M,[g])$ acts isometrically
for $\mu$, and, therefore, properly.

When $\mu$ does not define a distance, then Ferrand uses capacities to
construct conformally invariant functions on distinct quadruples or
triples of
points, no three of which are equal, depending whether $M$ is compact or not.  In the former case,
this function is essentially a log-cross-ratio, which can detect whether
distinct points are converging under a sequence of conformal transformations.
In the latter case, the function partly extends to the Alexandrov compactification $\hat{M} = M \cup \{
\infty \}$ and detects divergence to infinity, as well.
With these functions, she
reconstructs source-sink dynamics
 for unbounded sequences in
$\Conf(M,[g])$ when it acts nonproperly, and ultimately concludes $(M,[g]) \cong \mbox{\bf Euc}^n$ or ${\bf S}^n$.

\subsection*{PDEs proof, CR analogue}

When deforming a metric $g$ to one with desirable
properties, a natural choice is to restrict to conformal deformations.
The Yamabe Problem is a famous case: on any
compact Riemannian manifold $(M,g)$, there is $g' \in [g]$ with
constant scalar curvature.

In 1995, Schoen published a different proof of
Ferrand's theorem \cite{schoen.cr}, in the wake of his work on the completion of the
Yamabe Problem.   Given a Riemannian metric $g$ with scalar curvature
function $S_g$ and Laplace-Beltrami operator $\Delta_g$, if $fg \in [g]$ has constant scalar curvature $c$,
then $u = f^{\frac{n-2}{4}}$ satisfies
$$ - L_g u =
c \frac{n-2}{4(n-1)} u^{\frac{n+2}{n-2}}$$
where the conformal Laplacian
$$ L_g  : u \mapsto  \Delta_g u  - \frac{n-2}{4(n-1)} S_g u $$
is an elliptic operator satisfying a certain conformal invariance.
This equation can be applied to $f(x) = |\mbox{Jac}_x F|^{2/n}$ if $g$ has
constant scalar curvature and $F \in \Conf(M,[g])$.  Schoen's
arguments are based on the analytic properties of the elliptic
operator $L_g$.

\emph{CR}, or \emph{Cauchy-Riemann}, structures, model real hypersurfaces in complex
vector spaces.  On the unit sphere in $\BC^{n+1}$, for
  example, the tangent bundle carries a totally nonintegrable, or
  \emph{contact}, hyperplane distribution, equal at each $z \in
  S^{2n+1}$ to the maximal complex subspace of $T_z S^{2n+1} \subset
  {\bf C}^{n+1}$.  In general, a CR structure (nondegenerate, of
  hypersurface type) on an odd-dimensional manifold $M$ comprises a
  contact hyperplane distribution $\mathcal{D}$ and an almost-complex
  structure $J$ on $\mathcal{D}$, satisfying a
  compatibility condition: for any 1-form $\lambda$ on $M$ with
  $\mbox{ker } \lambda = \mathcal{D}$, 
the \emph{Levi form} on $\mathcal{D}$ given by $L_\lambda(u,v) := - d \lambda(u,v)$ is required
to obey $L_\lambda(Ju,Jv) = L_\lambda(u,v)$.  In this case, $L_\lambda$ is the
imaginary part of a Hermitian form on $\mathcal{D}$.  The real part of
this Hermitian form is a semi-Riemannian metric $B_\lambda$ on $\mathcal{D}.$

 A different choice of $\lambda$ gives
a metric conformal to $B_\lambda$, that is, related by a positive
function on $M$.  The CR structure is
\emph{strictly pseudoconvex} if any $B_\lambda$ is positive definite.
In this case, there is a subelliptic Laplace operator $L_\lambda$
on $C^\infty(M)$, sastifying a certain CR invariance.
Schoen proved, in analogy with his conformal results, that the
automorphisms of a strictly pseudoconvex CR manifold $M^{2n+1}$
act properly, unless
$M$ is CR-equivalent to $S^{2n+1} \subset {\bf C}^{n+1}$, or to
a Heisenberg group $\mbox{Heis}^{2n+1}$ carrying a certain
left-invariant CR structure.  S. Webster previously proved such a
theorem, in 1977, for $M$ compact with noncompact, connected
automorphism group, following some analogy with Obata's conformal proof.

\subsection*{Cartan connections proof, rank-one analogues}

The CR sphere can be identified as a conformal sub-Riemannian space with the visual boundary of complex hyperbolic
space ${\bf H}_{\BC}^{n+1}$.   The
quaternionic hyperbolic space ${\bf H}_{\bf H}^{n+1}$ and the Cayley hyperbolic
plane ${\bf H}^2_{\bf O}$ have visual boundaries diffeomorphic to spheres,
with respective differential-geometric structures invariant
by the isometries of the interior.  These boundaries of hyperbolic
spaces form a family;
they are parabolic homogeneous spaces ${\bf X} = G/P$ with $G$ simple of $\BR$-rank 1---see Table \ref{table.rank1}.
Cartan geometries modeled on one of these homogeneous spaces canonically
correspond to 
certain differential-geometric structures:

\begin{table}[h]
\begin{tabular}{|c|c|c|}
  \hline
  $ {\bf X}$ & $G$ & structure  \\
  \hline
  $\partial {\bf H}_{\BR}^{n+1}$ & $\mbox{PO}(1,n+1)$ & \small{conformal
                                                        Riemannian} \\
  $\partial {\bf H}_{\BC}^{n+1}$ & $\mbox{PU}(1,n+1)$ & \small{strictly
                                                    pseudoconvex CR} \\
  $\partial {\bf H}_{\bf H}^{n+1}$ & $\mbox{PSp}(1,n+1)$ & \small{contact
                                                         quaternionic} \\
  $ \partial {\bf H}^2_{\bf O}$ & $F_4^{-20}$ & \small{contact octonionic}  \\
                                        \hline
\end{tabular}
\caption{\label{table.rank1} rank-one parabolic geometries}
\end{table}
\noindent
Contact quaternionic and contact octonionic
structures were introduced by O. Biquard in 2000.

Let $M^n$  be
a manifold carrying one of the above structures, and suppose
that a group $H$ of transformations acts nonproperly---that is, there are
$h_k  \rightarrow \infty$ in $H$ and points $x_k \rightarrow x$ in $M$ such that $h_k.x_k
\rightarrow y$.  Let $\pi: \hat{M} \rightarrow M$ be the Cartan
bundle.
Choose lifts, $\hat{x}_k$ and $\hat{x}$ of $x_k$ and 
$x$, respectively, to $\hat{M}$, such that $\hat{x}_k \rightarrow \hat{x}$.

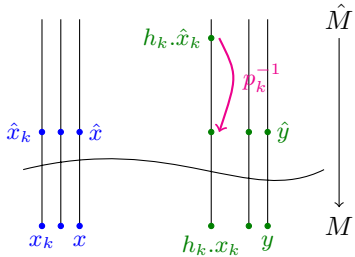
\begin{figure}
  \centering
\begin{tikzpicture}
\draw (0,-.5) .. controls (2,0) and (3,-1) .. (4,-.5)
(.25,1.5) -- (.25,-1.25)
(.5,1.5) -- (.5,-1.25)
(.75,1.5) -- (.75,-1.25) 
(2.5,1.5) -- (2.5,-1.25)
(3,1.5) -- (3,-1.25)
(3.25,1.5) -- (3.25,-1.25);
\filldraw[color=blue] (.25,-1.25) node[anchor=north]{\small $x_k$} circle (1pt);
\filldraw[color=blue] (.5,-1.25) circle (1pt);
\filldraw[color=blue] (.75,-1.25) node[anchor=north]{\small $x$} circle (1pt);

\filldraw[color=DarkGreen]
 (2.5,-1.25) node[anchor=north]{\small $h_k.x_k$} circle (1pt)
 (3,-1.25) circle (1pt)
 (3.25,-1.25) node[anchor=north]{\small $y$} circle (1pt);

\filldraw[color=blue] (.75,0) node[anchor=west]{\small $\hat{x}$} circle (1pt)
(.25,0) node[anchor=east]{\small $\hat{x}_k$} circle (1pt)
(.5,0) circle (1pt);
\filldraw[color=DarkGreen]
(2.5,0) circle (1pt)
(2.5,1.25) node[anchor=east]{\small $h_k.\hat{x}_k$} circle (1pt)
(3,0) circle (1pt)
(3.25,0) node[anchor=west] {\small $\hat{y}$} circle (1pt);
\draw[->,thick,color=magenta,shift={(3pt,0)}] (2.5,1.25) ..  node[anchor=west]{\small $p_k^{-1}$} controls (2.75,.75) .. (2.5,0); 
\draw[->] (4.2,1.25) node[anchor=south] {$\hat{M}$} -- (4.2,-1)
node[anchor=north] {$M$};
\end{tikzpicture}
\caption{\label{figure.holonomy} a holonomy sequence for $\{ h_k \}$ at $x$}
\end{figure}

Recall that the Cartan connection $\omega \in \Omega^1(\hat{M},\lieg)$
determines a parallelization of $T\hat{M}$, invariant by
the $H$-action on $\hat{M}$.  It follows that $H$ acts properly on
$\hat{M}$; in particular, $\{ h_k.\hat{x}_k \}$ diverges.  This
divergence is necessarily ``in the fiber direction''---that is, there
exist $p_k \rightarrow \infty$ in $P$ such that $\{ h_k.\hat{x}_k.p_k^{-1}
\}$ converges to some $\hat{y} \in \pi^{-1}(y)$.  Such a sequence
$\{p_k \}$ is called a \emph{holonomy sequence} for $\{ h_k \}$ at
$x$.  See Figure \ref{figure.holonomy}.  (For the initiated: Holonomy sequences are
related to $P$-valued cocycles for the $H$-action on $\hat{M}$, more
precisely, to \emph{placings} of $H$ into $P$.)

After passing to a subsequence, $\{ p_k \}$ has source-sink dynamics
on the model space ${\bf X}$.
Frances used the
Cartan connection to prove that $\{h_k \}$ has the same dynamical
behavior near $x$.  This led to his simultaneous
proof of Ferrand's and Schoen's theorems and of the extension to all
structures listed in Table \ref{table.rank1}: if the automorphism
group acts nonproperly, and $M$ is compact, then
it is geometrically isomorphic to ${\bf X}$.
If $M$ is noncompact, then it is isomorphic to the complement of a
point in ${\bf X}$.  This complement is
identified with a maximal connected, unipotent subgroup of $G$,
contained in $P$, via an orbit map.
For example, the strictly pseudoconvex CR structure on the Heisenberg group
appearing in Schoen's theorem can be defined this way.  Frances' proof
was a breakthrough, because his basic approach can be applied to a
very wide range of problems.

\subsection*{Pseudo-Riemannian analogue?}

The Ferrand--Obata Theorem is paradigmatic for the Zimmer--Gromov
Program.  In \cite{dag.rgs}, D'Ambra and Gromov ask whether
there is a 
pseudo-Riemannian analogue, in the case $M$ is compact.  Given that
pseudo-Riemannian isometry groups can act nonproperly, the suitable hypothesis
is that the conformal groups acts
essentially, meaning it does not preserve any metric in the
conformal class.  Alekseevsky
previously constructed an infinite family of non-conformally flat Lorentzian
metrics on ${\bf R}^n$ admitting an essential conformal flow.
The answer to D'Ambra and Gromov's question is ultimately ``no'';
nonetheless, we briefly describe here some positive results holding
under substantial assumptions either on the metric or on the
conformal group.

Let $(M^n,g)$ be a semi-Riemannian \emph{Einstein} space,
meaning the Ricci curvature satisfies
\begin{equation}
  \label{eqn.einstein}
  \mbox{Ric}_g = \frac{S_g}{n} g
  \end{equation}
In this case, the scalar curvature
$S_g$ is necessarily constant.  A notable precedent for the
Ferrand--Obata theorem was a result of Nagano--Yano for Riemannian
Einstein manifolds.  For simplicity, we state
the version for $M$ compact: if the conformal group contains a
nonisometric flow, then $(M,[g]) \cong ({\bf S}^n, [g_{+1}])$. 
For $g$ Einstein, of any signature, the
divergence $\sigma$ of a conformal vector field 
satisfies 
\begin{equation}
  \label{eqn.hessian.divergence}
  \nabla^2 \sigma = - \frac{S_g}{n(n-1)} \sigma g
\end{equation}
  When $M$ is compact, then $\sigma$ has extrema,
  which forces $S_g =0$ or $g$ definite; thus there is a
  parallel, isometric vector field, $\nabla \sigma$, or $(M,[g]) \cong ({\bf S}^n, [g_{+1}])$.
Local normal forms for metrics admitting nontrivial
solutions of the PDE (\ref{eqn.hessian.divergence}), due to Brinkmann,
rule out the first possibility and also give results for
noncompact, complete $(M,g)$, which were obtained by Kerckhove in 1988
and K\"uhnel--Rademacher in 2009.  Also in 2009, Kiosak--Matveev evaluated
(\ref{eqn.hessian.divergence}) along null geodesics for a clever
proof that any conformal rescaling of an indefinite $g$ which is also
Einstein must be constant, if $(M,g)$ is compact or complete.

For $(M^n,g)$ a compact, $(p,q)$-semi-Riemannian manifold and $H$ a 
simple group of conformal transformations, Zimmer proved in 1987 that $\mbox{rk}_{\BR} H \leq\mbox{min} \{p,q
\}+1$.  In light of his bound for isometry groups,
any $H$
attaining the conformal bound acts
essentially.  In this case, Bader--Nevo proved in 2002 that $H$ is locally
isomorphic to $\mbox{PO}(p+1,k)$, with $p+1 \leq k \leq q+1$.
Frances--Zeghib then strengthened their conclusion in 2005 to say
that $(M,g)$ is conformally equivalent to ${\bf
  S}^{p,q}$, up to double cover when $p,q \geq 2$, and up to $\BZ
\ltimes \BZ_2$-covers in the Lorentzian case.
These results support a pseudo-Riemannian Lichnerowicz Conjecture
when a large simple group acts.  Frances--Melnick proved an analogous
theorem when a large, connected, nilpotent group acts in 2010; here, large
means of maximal nilpotence degree.

In a recent advance, Pecastaing proved a similar result for a cocompact
lattice $\Gamma$ in a simple Lie group $H$, exploiting the proof of Zimmer's Conjecture.  The standard approach to a
$\Gamma$-action on a manifold $M$ is via the \emph{suspension} action
of $H$ on the diagonal quotient $H \times_\Gamma M$.
Brown--Fisher--Hurtado proved key relations on the Lyapunov
spectrum of this action in terms of the roots of $\lieh$.  When $\Gamma$ acts conformally on $M$, the
Lyapunov spectrum on the vertical distributions, tangent to the
$M$-fibers of the suspension, must be compatible with the spectrum of
the standard representation of $\mbox{CO}(p,q)$.
This strong restriction leads to the bound $\mbox{rk}_{\BR} H \leq
\mbox{min} \{ p,q \} +1$, and to uniform contracting behavior around a
point for a sequence in
$\Gamma$ when the bound is attained.  Pecastaing then proves that $(M,[g])$
is conformally flat; when $p,q \geq 2$, he obtains Frances--Zeghib's
conclusion that $(M,[g]) \cong {\bf S}^{p,q}$, up to double covers.

Finally, we come to the counterexamples: in
2012, Frances constructed non-conformally flat metrics on $S^1 \times
S^{n-1}$ of all signatures $(p,q)$ with $\mbox{min} \{p,q \} \geq 2$, admitting an
essential conformal flow.  We could still dream of a classification of
compact pseudo-Riemannian manifolds with essential conformal group,
but not of a Ferrand--Obata theorem in higher signature.  The CR story
has continued on a parallel course: Case--Curry--Matveev constructed compact CR manifolds
$(M,\mathcal{D}, J)$ with conformal metric $[B_\lambda]$ of arbitrary
higher---that is, non-Lorentzian---signature, which admit an essential
CR flow but are not CR-flat.  They also constructed a CR analogue of
Alekseevsky's noncompact, non-conformally flat, essential Lorentzian examples.

\subsection*{Lorentzian Lichnerowicz Conjecture}

We first describe some interesting examples showing that the global
conclusion of the Ferrand-Obata Theorem cannot hold in the
Lorentzian case.  Let $U = \mbox{Min}^{1,n-1} \backslash \{ 0 \}$.
The group $\Lambda = \{
2^k \mbox{Id}_n\}_{k \in \BZ}$ acts properly discontinuously and
cocompactly on $U$, commuting with $\mbox{CO}(1,n-1)$.  The compact
quotient $M = \Lambda \backslash U$ is called the conformal Lorentzian
\emph{Hopf manifold} and is diffeomorphic to $S^1 \times S^{n-1}$.
The conformal group is $S^1 \times
\mbox{O}(1,n-1)$, which means that $M$ is not conformally equivalent
to the double cover of ${\bf S}^{1,n-1}$.  It is, rather, a quotient of
the open, dense image of the conformal embedding of $U$ into ${\bf
  S}^{1,n-1}$, as in example 5.  This conformal action is also distinguished from the
M\"obius space by the special property that it preserves
an affine connection, descended from the standard affine connection
on $\BR^n$; the quotient connection is incomplete.

C. Frances constructed compact, Kleinian $(\mbox{PO}(2,n),{\bf
  S}^{1,n-1})$-spaces, for all $n \geq 3$, of infinitely many distinct
topological types, admitting an essential conformal flow.  They are quotients of the complement of a finite set of closed 
isotropic geodesics in ${\bf
  S}^{1,n-1}$ by the action of a Fuchsian group in
$\mbox{SL}_2(\BR)$ embedded in $\mbox{PO}(2,n)$.

The following Lorentzian Lichnerowicz Conjecture has, as yet, neither
been proved nor disproved:
\begin{conjecture}[LLC]
Let $(M^n,g)$ be a compact Lorentzian manifold, $n \geq 3$.  If
$\mbox{Conf}(M,[g])$ acts essentially, then $(M,g)$ is conformally flat.
  \end{conjecture}

The recent conformal D'Ambra Theorem of Melnick--Pecastaing verifies
the conjecture for $(M,g)$ real-analytic, with finite fundamental
group \cite{mp.confdambra}:

\begin{theorem}[Melnick--Pecastaing 2019]
For $(M^n,g)$ a compact, simply connected, real-analytic Lorentzian
manifold, for $n \geq 3$, $\Conf(M,[g])$ is compact.
    \end{theorem}
\noindent
A compact conformal group always preserves a metric in the conformal
class, the average of a given metric
in $[g]$ over the orbit.  Thus there are no essential conformal groups
for
$(M,g)$ as in the theorem, and the LLC is
vacuously true in this case.
Conversely, the LLC together with D'Ambra's
theorem would imply the above statement.
If $H = \Conf(M,[g])$ were
noncompact, D'Ambra's Theorem would say $H$ must be essential.  The LLC
would imply that $M$ is conformally flat. But there are no compact, simply connected, conformally flat Lorentzian
manifolds.  Indeed, for any manifold $M$ with a $(G,{\bf X})$-structure,
there is a \emph{developing map} from the
universal cover $\delta: \tilde{M} \rightarrow {\bf X}$.  If
$\tilde{M}$ is compact, then $\delta$ is a covering map, and the
$(G,{\bf X})$-structure is complete.  If $M = \tilde{M}$, this is not
possible with ${\bf X} = {\bf S}^{1,n-1}$, which has infinite
fundamental group.

\subsection*{Other conformal problems}

Two basic problems on compact pseudo-Riemannian manifolds, under the auspices of the VGC, are to classify the
homogeneous spaces with noncompact conformal group, and to classify
the possible connected noncompact conformal groups in general.

In the
Lorentzian case, the
LLC would indicate that the first problem comprises the analogous classification
of isometrically homogeneous Lorentzian spaces, and the classification
of compact, conformally flat Lorentzian manifolds with essential
conformal group.  The first problem was solved by Zeghib in \cite{zeghib.lorisom2}.  The second
problem may be approachable, but we do not currently know that all
such $(\mbox{PO}(2,n-1),{\bf S}^{1,n-1})$-spaces must be Kleinian. In
    higher signature, the first, isometric classification is not known. 

    For the second problem, we have the real-rank and nilpotence degree
    bounds mentioned above.   For compact Lorentzian manifolds, Pecastaing
    completed the classification of noncompact, semisimple, conformal
    groups.  They are necessarily of real-rank one or two; symplectic
    and exceptional groups do not occur.

With regard to conformal pseudo-Riemannian transformations, we are
still far from proving, but we are less far from starting the VGC
since D'Ambra and Gromov's influential manuscript of 30 years ago.
    
\section*{Projective Transformations}

Projective structures are geometrizations of systems of
second-order ODEs on $U \subset \BR^n$ of the form
\begin{equation}
  \label{eqn.proj.ode}
\frac{\ddot{\gamma}_1 + Q_1(\dot{\gamma})}{\dot{\gamma}_1} = \cdots =
\frac{\ddot{\gamma}_n + Q_n(\dot{\gamma})}{\dot{\gamma}_n}
\end{equation}
where $Q_i$ are quadratic forms on $\BR^n$. 
Note that the set of solutions is invariant under smooth reparametrization.
\'{E}. Cartan proved that for such a system, there are
torsion-free connections whose geodesics are the solutions.
 Two such connections $\nabla$ and $\nabla'$ have the same geodesics
 up to reparametrization.  They are related,
for vector fields $X$ and $Y$, by:
\begin{equation}
  \label{eqn.proj.equiv}
  \nabla'_X Y = \nabla_X Y + \nu(X) Y + \nu(Y) X
  \end{equation}
where $\nu$ is a 1-form.
A \emph{projective structure} on a manifold $M^n$ can thus be defined as an
equivalence class $[\nabla]$ of
torsion-free connections according to (\ref{eqn.proj.equiv}).  For $n \geq 2$, these correspond canonically
to Cartan
geometries modeled on ${\bf RP}^n$ and are rigid geometric structures
of order 2, formally quite similar to conformal structures.  The
projective group $\mbox{Proj}(M,[\nabla])$ comprises the
diffeomorphisms of $M$ preserving $[\nabla]$.

\subsection*{Strongly essential projective flows}

Consider the parabolic subgroup $P < \mbox{SL}_{n+1}(\BR)$ stabilizing
a line 
$\ell \subset \BR^{n+1}$.  It has Levi decomposition $P \cong \mbox{GL}_n(\BR)
\ltimes U$, where the unipotent radical $U$ acts
trivially on $\BR^{n+1}/\ell$, and thus with 
trivial differentials on $T_\ell {\bf RP}^n$.  An isomorphism
$U  \cong \BR^{n*}$ can be defined as follows: given $\xi \in \BR^{n*}
\cong (\BR^{n+1}/\ell)^*$, the element $u_\xi$ is determined by the
action on projective lines $\gamma_v$ through
$\ell$, as $v$ ranges over $\BR^{n+1}/\ell$; it effects the fractional linear reparametrization
\begin{equation}
  \label{eqn.proj.reparam}
  (u_\xi .\gamma_v)(s) = \gamma_v\left(\frac{s}{1 + \xi(v) s} \right)
  \end{equation}

Nagano--Ochiai considered a projective flow
$\{\varphi^t \}$ on an abritrary $(M, [\nabla])$ having a fixed point $x$
with $D_x\varphi^t$ trivial for all $t$.
For such a \emph{strongly essential} flow, there is $\xi \in T_x^*M$ such that each geodesic $\gamma_v$ is
reparametrized according to (\ref{eqn.proj.reparam}) with $t \xi$ in
place of $\xi$.  This can be calculated from the ODE
(\ref{eqn.proj.ode}) or can be seen via the Cartan
connection.
They then compute that, in a suitable framing 
along $\gamma_v(s)$,  the differential of
$\varphi^t$ is a scalar contraction by
$\lambda_s(t)  = (1+\xi(v) st)^{-1}$.  The flow by $\varphi^t$
contracts the projective Weyl curvature $W$ at $\gamma_v(s)$ by a factor of $\lambda_s(t)^2$. 
They conclude that $W_{\gamma_v(s)}$ vanishes for all
$v$ and all $s$.  A neighborhood of $x$ in $M$ is thus
projectively flat.  This 1986 paper was an early application of Cartan
connections for rigidity of transformation groups.  The technique
reappears and evolves in Frances' influential proof of the
Ferrand-Obata Theorem for rank-one parabolic geometries from above.

Local flatness results in the presence of strongly essential flows
have been proved for other parabolic geometries by Alekseevsky,
Frances, Frances--Melnick, \v{C}ap--Melnick, and Melnick--Neusser.

\subsection*{Riemannian Levi-Civita connections}

The projective analogue of the Lichnerowicz Conjecture would posit
that a compact affine manifold $(M^n,\nabla)$ with essential projective
group is equivalent to ${\bf RP}^n$, up to finite covers, with
\emph{essential} taken to mean not preserving a connection in $[\nabla]$.  Such a
statement has not yet received strong endorsement.

For $(M^n,\nabla)$ a Riemannian manifold with Levi-Civita
connection, however, the study of $
\mbox{Proj}(M,[\nabla])$ is classical---these are the transformations
preserving the set of unparametrized geodesics.
A Lichnerowicz-type conjecture appears in this context in the paper of
Nagano--Ochiai; in fact, they proved that a compact, Riemannian manifold admitting a
strongly essential projective flow is equivalent up to finite covers to ${\bf RP}^n$.

For metrics admitting essential
projective deformations, Levi-Civita found normal forms at generic
points.  Starting from these, Solodovnikov proved the following
theorem in 1969 for real-analytic Riemannian manifolds $(M^n,g)$, with $n \geq
3$.  The full theorem was proved by Matveev in 2007 \cite{matveev.riem.lich}.

\begin{theorem}[Matveev 2007]
 Let $H \leq \mbox{Proj}(M^n,[\nabla^g])$, for $M$ a compact, connected
 Riemannian manifold, $n \geq 2$, and $\nabla^g$ the Levi-Civita
 connection.  If $H$ is connected and does not act by affine
 transformations of $\nabla^g$, then $(M^n,g)$ is isometric to a sphere of constant
 curvature, up to finite covers.
 \end{theorem}
 \noindent
 The theorem in fact holds for $M$ noncompact and complete.  For $M$
 compact, Yano proved in 1952 that $\mbox{Aff}^0(M,\nabla^g) \cong
 \mbox{Isom}^0(M,g)$.  The hypothesis that $H$ is essential could thus
 be replaced with noncompactness.  We have another great
 examplar of the VGC.

 Matveev's proof uses techniques from integrable systems as well as
the crucial fact, presaged in work of Dini and Liouville, that the projective class of a Levi-Civita connection is a
finite-dimensional space. 
Any two Riemannian metrics $g$ and $g'$ are related by a field $L$ of
self-adjoint endomorphisms of $TM$.  The PDEs on $L$ corresponding to
$\nabla^g \sim \nabla^{g'}$ are not
linear.  However, a ``weighted'' version is: the metric $g'$ defined
for $u,v \in T_xM$ by
$$ g_x'(u,v) = \frac{1}{\det L_x} g_x(L^{-1}_xu,v)$$
is projectively equivalent to $g$ if and only if for all $u,v \in T_x
M$, for all $x \in M$,
  $$g_x((\nabla^g_u L)_xv,w) = \frac{1}{2} ( g_x(v,u) d\lambda_x(w) +
  g_x(w,u) d\lambda_x(v)) $$
where $\lambda = \mbox{tr } L$.  A \emph{BM-structure} is a smooth
field of self-adjoint endomorphisms satisfying this equation.  This
definition was made in 1979 by Sinjukov and was rediscovered by
Bolsinov--Matveev in 2003.


In Matveev's Theorem as stated, $H$ must be connected.  He gave a
counterexample on the 2-torus, based on a classfication due to Dini: Let $f$ be a
nonvanishing
function on $S^1$.  For
$$    g = \left( f(x) - \frac{1}{f(y)} \right) \left( \sqrt{f(x)} dx^2 
      +
      \frac{1}{\sqrt{f(y)}} dy^2 \right)$$
the involution exchanging $x$ and $y$ is an essential projective
transformation.
The reformulation of Matveev's Theorem with $H$ assumed noncompact,
however, holds even for $H$ disconnected.  Zeghib showed that 
$\left| \mbox{Proj}(M,[\nabla^g])/\mbox{Isom}(M,g) \right| \leq 2n$, unless $M$ is
isometric to a sphere,  in 2016; this bound was subsequently improved
to 2 by Matveev.

\subsection*{Pseudo-Riemannian Levi-Civita connections}

The definition of BM-structure applies as well to pseudo-Riemannian
metrics.  In this case, the
normal forms for the projectively related metrics are
more difficult to analyze.
Denote $\mathcal{L}_g$ the
space of BM-structures for a metric $g$.  Homothetic metrics are
affinely equivalent, so $\mbox{dim } \mathcal{L}_g \geq 1$ in general.

Let $(M^n,g)$ be a semi-Riemannian Einstein space, not necessarily
compact, with $n \geq 3$.   
Kiosak--Matveev
proved in 2009 that if $M$ is geodesically complete and admits essential
projective transformations then it is Riemannian and isometric to the sphere, up to finite covers.
They also proved that essential projective
transformations of $M$ give nonconstant solutions of the
following \emph{Gallot-Tanno equation}:
\begin{equation}
  \label{eqn.gallot.tanno}
  \begin{split}
  (\nabla^3 \sigma) (X,Y,Z) = c \left ( 2 (\nabla_x \sigma) \cdot
  g(Y,Z) \right. \\
 \left. + (\nabla_Y \sigma) \cdot g(X,Z) + (\nabla_Z \sigma) \cdot g(Y,X)
\right), \  c \in \BR
\end{split}
\end{equation}
where $\nabla = \nabla^g$; namely, $\sigma = \mbox{tr } L$ for $L$ the corresponding
BM-structure.  Here is a lovely analogy with the
divergence of conformal vector fields and (\ref{eqn.hessian.divergence}).

Now let $(M,g)$ be any compact, semi-Riemannian manifold.
In 2010, Matveev--Mounoud proved that
(\ref{eqn.gallot.tanno}) has nonconstant solutions only if
$(M,g)$ is Riemannian and isometric to a sphere, up to finite covers.
This fact was established for Riemannian metrics through work 
of Gallot, Tanno, and Hiramatu.  These two results combine to prove
the Projective Lichnerowicz Conjecture (PLC) for compact
$(M,g)$ when $g$ is Einstein. 

Kiosak--Matveev proved in 2010 that if $\mbox{dim } \mathcal{L}_g
\geq 3$, then there are
nonconstant solutions of (\ref{eqn.gallot.tanno}); this was
previously proved in the Riemannian case by Solodovnikov.  Thus the
PLC also holds when $(M,g)$ is compact and $\mbox{dim } \mathcal{L}_g
\geq 3$.

Bolsinov, Matveev, and Rosemann developed normal forms for
pseudo-Riemannian metrics admitting essential projective vector fields
that led to a proof in 2015 of the PLC for closed Lorentzian
manifolds.  The function $\sigma = \mbox{tr } L$ and, more generally, the nonconstant
eigenvalues of $L$, again play an important role.  The gradients of
the latter functions are coordinate vector fields in the normal forms.
The conclusion is ultimately reached by showing that $\sigma$ satisfies (5).
They simultaneously prove the Yano-Obata
Conjecture for \emph{c-projective structures} on closed manifolds of
dimension at least four.  This is the
K\"ahler version of metric projective structures and the PLC.
Lacking the space for this interesting topic, we refer
to the survey \cite{cemn.cproj},  and, for definite K\"ahler metrics, to papers of
Apostolov--Calderbank--Gauduchon.

\subsection*{The space $\mathcal{L}_g$}
Let $(M,g)$ be a compact pseudo-Riemannian manifold.  Assume $\mbox{dim }
\mathcal{L}_g = 2$.
The projective group of $[\nabla_g]$ acts linearly on $\mathcal{L}_g$
by a representation $\rho_g$.
In 2016 Zeghib showed $\mbox{Aff}(M,\nabla^g) =
\mbox{Isom}(M,g) = \ker \rho_g$, up to finite index subgroups.  The
quotient $\bar{H} = \mbox{Proj}(M,[\nabla^g])/\mbox{Isom}(M,g)$ acts
by fractional linear transformations on ${\bf P}(\mathcal{L}_g)$,
under a suitable identification with ${\bf RP}^1$.  He verifies
this action has finite kernel, and, by a beautiful analysis, proves that
$\bar{\rho}_g(\bar{H})$ lies in an elementary subgroup of
$\mbox{PGL}_2(\BR)$.  Therefore, up to finite index and finite
quotients,
$\bar{H}$ is isomorphic to a
subgroup of $\BR$.

Given any torsion-free connection $\nabla$, let $\mathcal{L}$ be the
space of semi-Riemannian metrics $g$ with $\nabla^g \sim \nabla$.  This is the space of solutions
of the \emph{metrizability problem} for $\nabla$, which has been
studied as early as the 1880s by Liouville.  Let $\mathcal{L}_0$ comprise
the symmetric 2-tensors $\eta \in \mbox{Sym}^2 TM$ for which $(\nabla
\eta)_0$, the trace-free component of $\nabla \eta$, vanishes.
Then $\mathcal{L}$ can be identified with the open set of nondegenerate elements of
$\mathcal{L}_0$.  If there is one $g \in \mathcal{L}$, then a
neighborhood of $g$ of dimension $\mbox{dim }\mathcal{L}_0$ is
contained in $\mathcal{L}$.

Eastwood--Matveev prolonged the equation $(\nabla \eta)_0$ $= 0$ and
found a nice characterization of the solutions. Let  $\hat{M}$ be the Cartan
bundle, and $P$, as usual, the stabilizer of a line $\ell$ in $\BR^{n+1}$.
The associated vector bundle $\mathcal{V}  = \hat{M} \times_P
\mbox{Sym}^2 \BR^{n+1}$ is called a \emph{tractor bundle}. 
The Cartan connection gives a projectively-invariant connection
$\nabla^\omega$ on $\mathcal{V}$.
Modifications of $\nabla^\omega$ by terms involving the
Cartan curvature give other meaningful, invariant
connections.  For one of these, given explicitly by Eastwood--Matveev,
the parallel sections are the solutions of the metrizability
equation, via a projection
corresponding to $\mbox{Sym}^2 (\BR^{n+1}/\ell)$.  Those satisfying a
nondegeneracy condition give metrics in the projective class.


The PLC for Levi-Civita connections in higher signature is believed to
be true, but has yet to be proved in the case $\mbox{dim }
\mathcal{L}_g = 2$.  For non-metric connections, a wide range of
questions on both projective and affine transformations are open.

\bibliography{karinsrefs}
\bibliographystyle{amsplain}







\end{document}